\documentclass[12pt,a4paper]{article}

\usepackage[english]{babel}
\usepackage{a4wide,epsfig,graphics,float,subfigure}
\usepackage{graphicx,color}
\usepackage{amsmath}
\usepackage{amssymb}
\usepackage{dsfont}
\usepackage{bbold}

\usepackage{algorithm2e,amsmath,lipsum,xcolor,caption}
\usepackage{datetime}
\usepackage{multirow}
\usepackage{bbm}
\newcommand{\mrm}[1]{\mathrm{#1}}
\usepackage{boxedminipage}
\usepackage{fancybox}       
\usepackage{fancyhdr}       
\newtheorem{theorem}{Theorem}[section]
\newtheorem{remark}{Remark}
    
    \newtheorem{proposition}[theorem]{Proposition}

\newenvironment{proof}{\begin{trivlist}
                       \item{\bf Proof:}}{\hfill $\square$
                     \end{trivlist}}

    \newcommand{\qed}{\nobreak \ifvmode \relax \else
          \ifdim\lastskip<1.5em \hskip-\lastskip
          \hskip1.5em plus0em minus0.5em \fi \nobreak
          \vrule height0.75em width0.5em depth0.25em\fi}

\begin{document}
 \thispagestyle{empty}

\vspace*{3cm}
\centerline{\huge A  regularization of  incompressible Stokes problem with }

\vspace{0.5cm}
\centerline{\huge Tresca friction condition}


%
%

\vspace*{1cm}
\begin{center}

 {\bf \sc A. Zafrar} 
\end{center}
\begin{center}

\noindent Department of Mathematics, Faculty of Sciences Dhar El Mahraz,B.P. 1796,
 Fez, Morocco

\noindent Email: .abderrahim.zafrar@usmba.ac.ma
%
\end{center}

\begin{abstract}

In the present article, we introduce and study a model addressing the Stokes problem with non-linear boundary conditions of the Tresca type.  We suggest a new procedure for regularizing incompressible fluid, i.e. we assume that the divergence $\nabla \cdot {\bf u}\in [-\epsilon,\,\epsilon]$ which leads to class of constrained  elliptic variational inequalities. We use a fixed point strategy to show the existence and uniqueness of a solution and we reformulate the problem as an equivalent constrained minimization problem.  An ADMM is applied to the minimization problem and some algorithm are provided.

\end{abstract}
\vspace{0.25cm}

%
\section{Introduction}
In applied mechanics, various applications require a certain dependence of the parameters in their resulting mathematical formulations (Lamé coefficients, friction coefficient,  dynamic and kinematic viscosity and many other parameters). It is well known (see for example \cite{ainsworth1997aspects,wihler2006locking,yi2022locking}) that the numerical approximation of such problems by  when low-order standard nodal-based displacement methods are used can fail to be robust when one or more parameters approach a some problem-dependent critical limit. This failure in the robustness of the finite element method is generally referred to as "locking" and appears in several forms. For example, very fine domains in shell and plate models can lead to so-called shear locking. In addition, the interaction between bending and membrane energies, which appears in shell theories, can lead to membrane locking. Finally, volume locking is observed in applications dealing with nearly incompressible materials.

In their attempt to overcome locking effects, a wide range of approaches have been developed. We mention here the mixed finite element methods and the nonconforming methods proposed in \cite{carstensen2000locking,guan2022locking}, and the higher-order methods in \cite{di2015hybrid}. All these methods have been quite extensively studied within an a priori context.
Another way of circumventing locking is the use of discontinuous Galerkin finite element methods \cite{yi2022locking}.  Stabilized finite element methods, e.g., of Galerkin/least-squares type is a different way to deal with the locking effect. This approach is similar to a mesh dependent relaxation of the incompressibility condition, as suggested \cite{brezzi1984stabilization}.

This work deals with the Stokes problem with a new regularization compressibility hypothesis.  This assumption expressed by the condition given by $\nabla {\bf u}\in [-\epsilon,\epsilon]$.  In other words, we allow to the divergence of velocity to vary in "small" interval that can be reduced to null value which corresponds to the incompressibility state.
It is well known that for problems approaching incompressible limits, standard finite element approaches suffer from suboptimal convergence due to a locking phenomena \cite{Migo2017}. 
Our assumption is very important and more realistic since some elastic some fluids are slightly compressible and this can be linked to the locking effect as we will discover later.

According to  \cite{DemSuq}, materials
are classified into two categories: soft and hard. Concerning the first category ( aluminium for example), the stress
increment required to produce a specified strain increment diminishes
with increasing deformation. For the second one (cellular rubber for example),  the stress increment grows as the deformation proceeds.  An extreme case of a soft material is the rigid-plastic one.  We consider an elastic ideally locking fluids and then the constitutive law may be modelled by the (convex) subdifferential of the indicator function of a convex set $\mathcal T$ which characterizes the locking constraints.  This kind of model can be viewed as class of fluids  with  limited compressibility introduced  by Prager (see for instance \cite{DemSuq} and the references therein).  For such  class  of   fluids  the locking constraint  acts  only on the volumic part of  the  strain: 

\begin{equation}
\mathcal{T}=\left\lbrace \varepsilon : tr(\varepsilon )\leq \epsilon\right\rbrace 
\end{equation}
where $tr(\varepsilon )$ denotes the trace of the symmetric tensor $\varepsilon=(\varepsilon_{ij})$ and we have
\begin{equation}
\varepsilon ({\bf u}) \in \mathcal{T}\Rightarrow -\epsilon \leq \nabla \cdot {\bf u}\leq \epsilon \label{impl}
\end{equation}
We note that incompressible elastic materials  are included in this case for $\epsilon=0$.

On the other hand,
no-slip hypothesis at fluid-wall interface is assumed which is showed to be a good agreement with experimental observations for both newtonian and non-newtonian fluid \cite{magnin87}.  
No-slip hypothesis means that the velocity at the wall is not zero.
This boundary condition was firstly adopted in a numerical simulation of a flow by Doltsini {\it et al.} \cite{Doltsini87} and Fortin \cite{fortin91}. Since that, many papers were published simulating various flows with such boundary conditions (see \cite{rao99} and references therein).  More recently, based on the penality method, an error estimates for the Stokes problem with Tresca boundary conditions  are obtained \cite{kaitai08}. 

Our aim in this paper is to contribute to the mathematical modelling and numerical analysis of incompressible Stokes problem with new kind of regularisation and with Tresca boundary conditions.  We first perform the existence and uniqueness of the solution to this problem.  To do this,  we consider a convergent fixed point scheme.  By using this approach, 
Our second purpose is to derive a minimization problem that can be simple to handle with standard numerical method.  

The paper is organized as follows. In the next Section, we introduce the constitutive  equations modelling the Stokes problem with nearly incompressibility condition. Section 3 is devoted to the existence and uniqueness of the solution to the regularized problem. Moreover, we reformulate an equivalent minimization problem. Section 4 and 5 are devoted the numerical analysis,  we apply an ADMM to split the problem and makes the constraints explicit,   and we provide some algorithms.
\section{Nearly incompressible Stokes problem with Tresca friction condition}
The aim of this section is to introduce some notations and recall some functional structures that are necessary for our study. We will consider an open bounded set $\Omega\subset \mathbb{R}^d$, $d\ge2$,  with regular boundary $\partial \Omega:=\Gamma$ which is the divided into  two disjoint parts $\Gamma_d $ and $\Gamma_f $ (eventually empty). 

The space $L^2(\Omega)$ of square integrable functions in $\Omega$  is equipped with the norm:
$$
\displaystyle \Vert p\Vert^2_0=  \int_\Omega|p|^2dx\, ,\qquad \forall p\in L^2(\Omega),
$$
while
 $L^2_0(\Omega)$ stands for the closed subspace of $L^2(\Omega)$ defined by:
$$
\displaystyle L^2_0=\left\lbrace  p\in  L^2(\Omega) \mbox{ such that }  \int_\Omega p\,dx\, =0\right\rbrace.
$$

\noindent The commonly use  notation $H^m(\Omega)$, $m\ge1$  stands for the standard Sobolev space, endowed with the usual norm:
$$
\displaystyle||{\bf u}||_m^2	=\sum\limits_{0\leq |\alpha| \leq m}||\partial^\alpha {\bf u}||^2_0,
$$
where $\alpha$ is a multi-index. 

%

Let $H^{-1}(\Omega)^d$ the dual of $H^1(\Omega)^d$ and $V\subset H^1(\Omega)^d$  be a subspace of functions vanishing on the open portion $\Gamma_d$, with non vanishing measure, i.e.:
$$
V=\left\lbrace {\bf v}\in H^1(\Omega)^d \mbox{ such that } {\bf v}_{|_{\Gamma_d}}=0\right\rbrace .
$$
%
%

We consider the following incompressible Stokes problems with nonlinear boundary condition of Tresca friction type:
\begin{eqnarray}
\label{Stok_tres_cont0}
\left\|
\begin{array}{cc}
-\nabla\cdot (\nu\varepsilon( {\bf u})) + \nabla p	=\textbf{f}			& \mbox{ in }\Omega\\
\nabla \cdot {\bf u}					= {\bf 0}				& \text{ in }\Omega \\
{\bf u}						=\textbf{0}			& \mbox{on }\Gamma_d\\
{\bf u}_{\vec{n}}=0\,	\text{ and }
\begin{cases}|\sigma_\tau|< \xi \Rightarrow  \quad {\bf u}_\tau	=\textbf{0},		\\
|\sigma_\tau|= \xi \Rightarrow \exists \kappa>0, \, {\bf u}_\tau	=-\kappa \sigma_\tau,			
\end{cases}
&\mbox{on }\Gamma_f,

\end{array}
\right.
\end{eqnarray}

 and the new regularized incompressible Stokes problem:
\begin{eqnarray}
\label{Stok_tres_cont}
\left\|
\begin{array}{cc}
-\nabla \cdot(\nu\varepsilon( {\bf u})) + \nabla p	=\textbf{f}			& \mbox{in }\Omega\\
\nabla \cdot {\bf u}					\in [-\epsilon, \epsilon] 				&  \\
{\bf u}						=\textbf{0}			& \mbox{on }\Gamma_d\\
{\bf u}_{\vec{n}}=0\,	\text{ and }

\begin{cases}
|\sigma_\tau|< \xi \Rightarrow  \quad {\bf u}_\tau	 =\textbf{0},			\\
|\sigma_\tau|= \xi \Rightarrow \exists \kappa>0, \, {\bf u}_\tau	=-\kappa \sigma_\tau,
\end{cases}
			&\mbox{on }\Gamma_f,
\end{array}
\right.
\end{eqnarray}

 where $\kappa$ is the friction coefficient.  $\Gamma_d $ is subjected to no-slip boundary condition (Dirichlet boundary condition) while $\Gamma_f$ is where the fluid may slip, $\nu>0$ is the kinematic viscosity  and   $\xi$ is a non-negative function in $L^2(\Gamma)$.  The linearized strain tensor is denoted by $\displaystyle \varepsilon({\bf u})=\frac{1}{2}(\nabla {\bf u}+(\nabla {\bf u})^\top)$, where $(\cdot)^\top$ is the transpose operator. 
We denote by $\vec{n}$ the outward unit normal to $\Gamma$ and ${\bf u}_{\vec{n}}$, respectively ${\bf u}_\tau$, the normal and the tangential, component of ${\bf u}$.

Let us define the Following spaces
$$
{ \bf V}=\{ {\bf v}\in V ,\\ {\bf v}_{\vec{n}}=0\text{ on }\Gamma_f\}, \qquad\quad {\bf H}_{\mrm{div}}=\{ {\bf v}\in {\bf V}\ ,\; \nabla \cdot{\bf v}=0 \mbox{ in }\Omega\},
$$
$${\bf V}_{\mrm{div}}=\{ {\bf v}\in {\bf V}\ ,\;\left( \nabla \cdot{\bf v},q\right)=0\;\;\forall q\in L_0^2(\Omega)\},$$

$$ {\bf H}_{\mrm{div},\epsilon}=\{ {\bf v}\in {\bf V}\ ,\; \nabla \cdot {\bf v}\in [-\epsilon,\;\epsilon] \mbox{ a.e.  in }\Omega\},$$
the following bilinear forms
$$\displaystyle \mathcal{A}({\bf u},{\bf v})= \int_{\Omega}\nu \varepsilon({\bf u})\colon \varepsilon({\bf v})dx\, ,\qquad \mathcal{B}(q,{\bf v})=\displaystyle \int_\Omega q \,\nabla\cdot{\bf v}dx,$$  
and the friction map
$$\displaystyle j({\bf v})=\int_{\Gamma_d}\xi|{\bf v}_\tau|ds \, \qquad \forall{\bf u}, {\bf v} \in {\bf V}.$$

Let's now handle the constraint assumed for the divergence of the velocity ${\bf u}$ that  appears in \eqref{Stok_tres_cont}. To do so,  we give an equivalent set to ${\bf H}_{\mrm{div},\epsilon}$.  We are seeking  a map $\mathcal{L}(\cdot)$ and then for
$$\mathbb{L}=\left\lbrace \varepsilon \in \mathbb{S}^d, \quad \mathcal{L}(\varepsilon)\leq 0\right\rbrace$$
such that ${\bf H}_{\mrm{div},\epsilon}$ is reformulated is as follows:
$${\bf H}_{\mrm{div},\epsilon}=\left\lbrace {\bf v}\in  {\bf V},\quad \varepsilon({\bf v})\in \mathbb{L},\;\; \mbox{ a.e.  in }\Omega\right\rbrace$$
where $\mathbb{S}^d$ is the usual space of symmetric second order tensors. Following \cite{Migo2017}, a suitable map $\mathcal{L}(\cdot)$ is defined as follows: 
\begin{equation}
\mathcal{L}(\varepsilon)=tr(\varepsilon)-\epsilon
\end{equation}
where $tr(\cdot)$ stands for the trace of a given tensor. 

It is straightforward that 
\begin{equation}\label{intersection}
\displaystyle {\bf H}_{\mrm{div}} =\bigcap_{\epsilon\searrow 0}{\bf H}_{\mrm{div},\epsilon} 
\end{equation}
 and for all $\epsilon>0$, $ {\bf H}_{\mrm{div},\epsilon}$ are closed and convex and $ \mathbf{0} \in  {\bf H}_{\mrm{div},\epsilon}$.


We recall the following result  (see \cite{saidi04}) to derive the variational problem.
%
\begin{proposition}\label{prop1}
\begin{equation}
\label{equiv_tresca}
\left\{
\begin{array}{rcllrr}
|\sigma_\tau|<\xi \, \Rightarrow {\bf u}_\tau					&=&0\\
|\sigma_\tau|=\xi \, \Rightarrow {\bf u}_\tau					&=&-\kappa \sigma_\tau \quad \mbox{ on }\Gamma_f
\end{array}
\right.
\qquad \Longleftrightarrow \qquad 
\sigma_\tau.{\bf u}_\tau +\xi|{\bf u}_\tau|=0 \mbox{ on }\Gamma_f
\end{equation}
\end{proposition}
%

In this paper we consider  the constitutive law where the stress tensor is  defined in the form
$$
\sigma\in 2\nu\varepsilon({\bf u})-p \mathbbm{I}_d+\partial\mathbbm{1}_{\mathbb{L}}(\varepsilon ({\bf u}))
$$
where  $p$ is the pressure, $\mathbbm I_d$ stands for the identity tensor, $\nu>0$ represents the kinematic fluid viscosity.

On the other hand,  $\sigma\in 2\nu\varepsilon({\bf u})-p \mathbbm{I}_d+\partial\mathbbm{1}_{\mathbb{L}}(\varepsilon ({\bf u}))$ implies that there exists some  selection of $\varepsilon({\bf u})$ in $\partial \mathbb{L}$ denoted by $\ell(\varepsilon({\bf u}))$ such that 
$$\sigma= 2\nu\varepsilon({\bf u})-p \mathbbm{I}_d+\ell(\varepsilon ({\bf u})).$$
Next, from the first equation and the second in \eqref{Stok_tres_cont} and taking into account \eqref{impl} we derive the following equation
\begin{equation}
\int_\Omega \nu\varepsilon({\bf u})\varepsilon({\bf v})dx-\int_\Omega  p\, \nabla\cdot {\bf v}dx+\int_\Omega \ell(\varepsilon ({\bf u}))\varepsilon({\bf v})dx=\int_\Omega {\bf f v}dx +\int_{\Gamma}\sigma {\vec{n}} {\bf v}ds,\qquad \forall {\bf v} \in {\bf H}_{\mrm{div},\epsilon},\,\forall \epsilon>0\label{vareq}
\end{equation}
for some $\ell(\varepsilon (u))\in \partial \mathbbm{1}_{\mathbb{L}}(\varepsilon ({\bf u}))  $. 

Now,  we deal with the inclusion boundary conditions. Recall that $\sigma= 2\nu\varepsilon({\bf u})-p \mathbbm{I}_d+\ell(\varepsilon ({\bf u}))$ for some $\ell(\varepsilon ({\bf u}))\in \partial \mathbbm{1}_{\mathbb{L}}(\varepsilon ({\bf u}))  $, then 
\begin{equation}\label{xiest}
 0\geq \mathbbm{1}_{\mathbb{L}}(\varepsilon ({\bf v}))- \mathbbm{1}_{\mathbb{L}}(\varepsilon ({\bf u}))\geq \ell(\varepsilon ({\bf u}))(\varepsilon ({\bf v})-\varepsilon ({\bf u})).
\end{equation}

 Thus an equivalent formulation of the above  problem is:  
\begin{eqnarray}
\label{ineq_var_stokes_gen_sigm1}
\left\|
\begin{array}{l}
\mbox{For all }\epsilon >0, \mbox{ find } ({\bf u},p) \in {\bf H}_{\mrm{div},\epsilon}\times L^2_0(\Omega)\mbox{ such that :}\, \\\\\
\displaystyle \mathcal{A}({\bf u},{\bf v}-{\bf u})+\mathcal{B}(p, {\bf v}-{\bf u})+j({\bf v})-j({\bf u})\geq F({\bf v}-{\bf u}),\forall \, {\bf v}\in {\bf H}_{\mrm{div},\epsilon}.
\end{array}
\right.
\end{eqnarray}
where $\displaystyle F({\bf v})=\int_\Omega {\bf f v}dx +\int_{\Gamma}\sigma {\vec{n}} {\bf v}ds$.

\section{Existence results}
The variational problem \eqref{ineq_var_stokes_gen_sigm1} maybe considered as a control problem  with suitable cost function, where the pressure $p$ is the control. Based  on this idea,  to prove the existence of of solution to the this problem, we will consider a fixed point strategy.

 Since $p$ maybe considered as a Lagrange multiplier,  we consider the following "convergent" iterative scheme \eqref{itersch}, where we fix the pressure in the  inequality and we give an update of it.
\begin{eqnarray}
\label{itersch}
\left\|
\begin{array}{l}
\mbox{For all }\epsilon >0, \mbox{ find } ({\bf u}^{n+1},p^{n+1}) \in {\bf H}_{\mrm{div},\epsilon}\times L^2_0(\Omega)\mbox{ such that :}\, \\\\\
\displaystyle \mathcal{A}({\bf u}^{n+1},{\bf v}-{\bf u}^{n+1})+j({\bf v})-j({\bf u}^{n+1})\geq F({\bf v}-{\bf u}^{n+1})+\int_\Omega  p^{n}\, \nabla \cdot({\bf v}-{\bf u}^{n+1})dx,\forall \, {\bf v}\in {\bf H}_{\mrm{div},\epsilon},\\
p^{n+1}=p^{n}+\Pi_{L^2}(\nabla\cdot ({\bf u}^{n+1})),
\end{array}
\right.
\end{eqnarray}
where $\Pi_{L^2}$ is the orthogonal projection on $L_0^2(\Omega)$.
   
The following result gives the convergence of the iterative scheme \eqref{itersch}. This fixed point strategy shows that the problem \eqref{ineq_var_stokes_gen_sigm1} admits a unique solution, also it allows us to consider the problem as quasi-variational inequality.

\begin{proposition}
The problem \eqref{ineq_var_stokes_gen_sigm1} and \eqref{itersch} admit a unique solution and the solution $({\bf u}^{n},p^n)$ provided by this iterative scheme  converges strongly in ${\bf H}_{\mrm{div},\epsilon}\times L^2_0(\Omega)$  to the solution to the problem \eqref{Stok_tres_cont}.
\end{proposition}
\begin{proof}
The existence and uniqueness of the solution to the iterative scheme \eqref{itersch} is straightforward. 
In fact,  for fixed $p^{n}\in L^2(\Omega)$ such that $\displaystyle\int_{\Omega} p^ndx =0$, the bilinear form form $\mathcal{A}(\cdot,\cdot)$ is continuous and coercive, the linear form ${\bf v}\longmapsto F({\bf v})+\displaystyle\int_\Omega p^n\nabla \cdot {\bf v}dx$ is bounded and $j(\cdot)$ is proper, lower semi-continuous and strictly convex. Then,by standard result in variational and quasivariational inequalities, the solution to \eqref{itersch} exits and it is unique. Let us show the strong convergence of $({\bf u}^n,p^n)_n$ in ${\bf H}_{\mrm{div},\epsilon}\times L^2_0(\Omega)$.

 In the one hand, given $p^n$ and let us take ${\bf v}={\bf 0}$ in the variational inequality, and because of $j({\bf v})\geq 0$ for all ${\bf v}$,  we get:
\begin{eqnarray}
\mathcal{A}({\bf u}^{n+1},{\bf u}^{n+1})+j({\bf u}^{n+1})\leq F({\bf u}^{n+1})+\int_\Omega  p^{n}\, \nabla \cdot{\bf u}^{n+1}dx
\end{eqnarray}
implies that 

\begin{eqnarray}
\mathcal{A}({\bf u}^{n+1},{\bf u}^{n+1})\leq F({\bf u}^{n+1})+\int_\Omega  p^{n}\, \nabla \cdot{\bf u}^{n+1}dx
\end{eqnarray}
since  $\mathcal{A}(\cdot,\cdot)$ is coercive, there exists a positive constant $\alpha$ such that 
\begin{equation}
\alpha \Vert {\bf u}^{n+1}\Vert^2 \leq \mathcal{A}({\bf u}^{n+1},{\bf u}^{n+1}),
\end{equation}
this implies that
\begin{eqnarray}
\Vert {\bf u}^{n+1}\Vert_1^2 \leq \dfrac{1}{\alpha}\left(\Vert F\Vert_0+C\Vert  p^{n}\Vert_0\right)\Vert {\bf u}^{n+1}\Vert_1
\end{eqnarray}
this shows that the sequence $({\bf u}^{n+1})_n$ is bounded and we can extract a weakly convergence subsequence denoted  by the same notation. On the other hand,
we have 

\begin{equation}
p^{n+1}=p^{n}+\Pi_{L^2}(\nabla \cdot{\bf u}^{n+1}).
\end{equation}
and since the projection operator $\Pi_{L^2}$ is self-adjoint and the adjoint of the divergence $\nabla\cdot(\cdot)$ is $-\nabla(\cdot)$,  we get

\begin{equation*}
\left( p^{n+1}, q\right) =\left( p^{n}, q\right) -\left({\bf u}^{n+1}, \nabla \Pi _{L^2}(q)\right)\qquad \forall q\in L^2_0(\Omega),
\end{equation*}
by passing to limit, we obtain
\begin{equation}
\lim_{n\rightarrow+\infty} \left( p^{n+1}, q\right) =\lim_{n\rightarrow+\infty} \left( p^{n}, q\right) -\left({\bf u}^*, \nabla \Pi_{L^2} (q)\right)\qquad \forall q\in L^2_0(\Omega),
\end{equation}
where ${\bf u}^*\in {\bf H}_{\mrm{div},\epsilon}$ is the weak limit of ${\bf u}^{n+1}$.  
\\
Now, for $\epsilon=\dfrac{1}{n}$, we get: 
\begin{equation}
\lim_{n\rightarrow+\infty} \left( p^{n+1}, q\right) -\lim_{n\rightarrow+\infty} \left( p^{n}, q\right)=0, \qquad \forall q\in L^2_0(\Omega).
\end{equation}
since ${\bf u}^*\in {\bf H}_{\mrm{div},\frac{1}{n}}$.

And this implies that 
\begin{equation}
\lim_{n\rightarrow+\infty}( p^{n+1} - p^{n})=0,
\end{equation}
as consequence, the sequence $(p^n)_n$ converges strongly in $L^2_0(\Omega)$ to $p^*$  the strong limit,  since $L^2(\Omega)$ is a Banach space. 

By the pseudo-monotony of $\mathcal{A}(\cdot,\cdot)$,  the $H^1$-Sobolev embedding (Rellich–Kondrachov theorem see \cite{adams2003sobolev}),   the strong continuity of the divergence and the gradient operators and since $j(\cdot )$ is weakly semi-continuous,  $({\bf u}^*,p^*)$ is a solution to the problem  \eqref{ineq_var_stokes_gen_sigm1} by passing to the limit in \eqref{itersch} where $n\longrightarrow +\infty$.

Let us show the strong convergence of the sequence $({\bf u}_n)_n$ which implies the uniqueness of the solution. Let $({\bf u},p)$ be a  solution to the problem: 

\begin{eqnarray}
\label{stkoesnear}
\left\|
\begin{array}{l}
\mbox{For all }\epsilon >0, \mbox{ find } ({\bf u},p) \in {\bf H}_{\mrm{div},\epsilon}\times L^2_0(\Omega)\mbox{ such that:}\, \\\\\
\displaystyle \mathcal{A}({\bf u},{\bf v}-{\bf u})+j({\bf v})-j({\bf u})\geq F({\bf v}-{\bf u})+\int_\Omega  p\, \nabla\cdot ({\bf v}-{\bf u})dx,\forall \, {\bf v}\in {\bf H}_{\mrm{div},\epsilon}.
\end{array}
\right.
\end{eqnarray}
By taking ${\bf v}={\bf u}^{n+1}$ in the quasivariational problem \eqref{stkoesnear} and ${\bf v}={\bf u}$ in \eqref{itersch} and adding the resulting inequalities we obtain:
\begin{eqnarray}
\mathcal{A}({\bf u}-{\bf u}^{n+1},{\bf u}-{\bf u}^{n+1})&\leq & \int_\Omega (p-p^n)(\nabla \cdot({\bf u}-{\bf u}^{n+1}))dx\nonumber\\
&&\text{ by Cauchy-Schwarz inequality}\nonumber\\
&\leq &C \Vert \nabla \cdot({\bf u}-{\bf u}^{n+1})\Vert_0 \Vert p-p^n\Vert_0\nonumber \\
&\leq &C \Vert \nabla({\bf u}-{\bf u}^{n+1})\Vert_0 \Vert p-p^n\Vert_0\nonumber \\
&\leq &C \Vert {\bf u}-{\bf u}^{n+1}\Vert_1 \Vert p-p^n\Vert_0 
\end{eqnarray}
where $C>0$ is some positive constant. The bilinear form $\mathcal{A}(\cdot,\cdot)$ is coercive,  i.e.  there exists $\alpha>0$ such that 
\begin{equation}
\alpha \Vert{\bf u}-{\bf u}^{n+1} \Vert_1 ^2\leq \mathcal{A}({\bf u}-{\bf u}^{n+1},{\bf u}-{\bf u}^{n+1})
\end{equation}
 this implies that:
\begin{equation}
\alpha \Vert{\bf u}-{\bf u}^{n+1} \Vert ^2_1\leq C \Vert {\bf u}-{\bf u}^{n+1}\Vert_1 \Vert p-p^n\Vert _0
\end{equation}
and hence
\begin{eqnarray}
\Vert {\bf u}-{\bf u}^{n+1}\Vert_1 
&\leq &(C/\alpha)  \Vert p-p^n\Vert _0.
\end{eqnarray}

This last inequality provide the strong convergence of the sequence $({\bf u}^{n})_n$, solution to the fixed point iterative scheme \eqref{itersch} to the solution the nearly incompressible problem \eqref{Stok_tres_cont}.
\end{proof}
\begin{remark}
 In the proof of the above result,  for the choice of $\epsilon=\frac{1}{n}$, the  strong limit of the sequence $({\bf u}^n,p^n)_n$ is the solution of the problem \eqref{Stok_tres_cont0}.  In other words, we have shown the existence of the solution of \eqref{ineq_var_stokes_gen_sigm1}, its convergence to the solution to the solution of \eqref{Stok_tres_cont0} and we have provided a strong tool to provide the solution numerically.
 \end{remark}
 \begin{remark}
With fixed pressure $p^n$, is clear that the velocity is the solution of quasivariational inequality arising in elliptic partial differential equations. In the next section, we will see that this iterative scheme combined with the ADMM  will perfectly address the locking effect.  More precisely, we will show that at each iteration, the velocity is the solution to a linear elastic problem with bilateral contact and Tresca friction. 

\end{remark}

Let us denote  $\displaystyle F_{p^n}({\bf v})=F({\bf v})+\int_\Omega  p^{n}\, \nabla\cdot {\bf v}dx$,  we simply express the problem \eqref{itersch} in the following form:
\begin{eqnarray}
\label{itersch1}
\left\|
\begin{array}{l}
\mbox{For all }\epsilon >0, \mbox{ find } ({\bf u}^{n+1},p^{n+1}) \in {\bf H}_{\mrm{div},\epsilon}\times L^2_0(\Omega)\mbox{ such that :}\, \\\\\
\displaystyle \mathcal{A}({\bf u}^{n+1},{\bf v}-{\bf u}^{n+1})+j({\bf v})-j({\bf u}^{n+1})\geq F_{p^n}({\bf v}-{\bf u}^{n+1}),\forall \, {\bf v}\in {\bf H}_{\mrm{div},\epsilon}.\\
p^{n+1}=p^{n}+\Pi_{L^2}(\nabla\cdot {\bf u}^{n+1}).
\end{array}
\right.
\end{eqnarray}

%

We can reformulate the above symmetric problem  \eqref{itersch1} as the  following minimization problem: 

For all $\epsilon >0$,  find $ {\bf u}\in 	\mathcal{A}_{ad,\epsilon}:={\bf H}_{\mrm{div},\varepsilon}$   solution to the minimization problem:
\begin{eqnarray}
\label{cst_minimi1}
\begin{array}{ll}

\displaystyle\inf_{{\bf v}\in \mathcal{A}_{ad,\epsilon} } J({\bf v})+ j({\bf v}), 
\end{array}
\end{eqnarray}
where $J({\bf u}):=\frac{1}{2}\mathcal{A}({\bf u},{\bf u}))-F_{p^n}({\bf u})$, $\mathcal{A}_{ad,\epsilon}$ is the admissible set and ${\bf u}^{n+1}$ is denoted simply by ${\bf u}$.  

The existence and uniqueness of the solution to the problem \eqref{cst_minimi1} follows from the fact that  the bilinear form $a(\cdot,\cdot)$ is continuous and coercive, the function $j(\cdot)$ is proper, strictly convex and lower semi-continuous and the linear form $F_{p^n}(\cdot)$ is bounded. (see \cite{ekel74}).  

\section{ADMM splitting }
The purpose of this section is to provide an ADDM splitting procedure.  To do this, we give  saddle point formulation of the problem \eqref{cst_minimi1}.   So,  since $j(\cdot)$ is not differentiable,  let us introduce an auxiliary unknown $\Phi$ and reformulate the problem \eqref{cst_minimi1} as follows:

\begin{eqnarray}
\label{cst_minimi3}
\left\|
\begin{array}{l}
\forall\epsilon >0, \mbox{ find } ({\bf u},\Phi)\in 	{\bf H}_{\mrm{div},\varepsilon}\times L^2(\Gamma_f) \mbox{ such that :}\\\\
\displaystyle J({\bf u})+ j(\Phi) \leq J({\bf v})+ j(\Psi), \quad \forall \, ({\bf v},\Psi)\in {\bf H}_{\mrm{div},\varepsilon}\times L^2(\Gamma_f) ,\\\\
{\bf u}_\tau-\Phi=0\quad\mbox{ on } \Gamma_f.
\end{array}
\right.
\end{eqnarray}

We need to make the constraint  on the divergence explicit, so the minimization problem \eqref{cst_minimi3} can be written in following equivalent form:
\begin{eqnarray}
\label{cst_minimi2}
\left\|
\begin{array}{l}
\mbox{For all }\epsilon >0, \mbox{ find } ({\bf u},\Phi_\tau)\in 	{\bf V}\times L^2(\Gamma_f)  \mbox{ such that :}\\\\
\displaystyle J({\bf u}) + j(\Phi_\tau) +
\mathds{1}_{[-\epsilon,\epsilon]}(\nabla \cdot {\bf u})
\leq J({\bf v}) + j(\Psi_\tau)+\mathds{1}_{[-\epsilon,\epsilon]}(\nabla \cdot {\bf v}) ,\qquad \forall \, ({\bf v},\Psi_\tau)\in {\bf V} \times L^2(\Gamma_f) ,\\
\mbox{subject to: } 
{\bf u}_\tau-\Phi_\tau=0\quad\mbox{ on } \Gamma_f,
\end{array}
\right.
\end{eqnarray}
where $\mathds{1}_{[-\epsilon,\epsilon]}(\cdot)$ is the indicator function.  

We introduce a new auxiliary in order to handle the constraint induced by the indicator function, i.e.:
\begin{eqnarray}
\label{cst_minimi2}
\left\|
\begin{array}{l}
\forall\epsilon >0, \mbox{ find } ({\bf u},\Lambda,\Phi_\tau)  \mbox{ such that :}\\\\
\displaystyle J({\bf u}) + j(\Phi_\tau) +
\mathds{1}_{[-\epsilon,\epsilon]}(\Lambda)
\leq J({\bf v})+ j(\Psi_\tau)+\mathds{1}_{[-\epsilon,\epsilon]}(\phi) ,\qquad
\forall \, ({\bf v},\phi,\Psi_\tau)\\\\
{\bf u}_\tau-\Phi_\tau=0\quad\mbox{ on } \Gamma_f,\\
\nabla \cdot {\bf u}-\Lambda=0\quad\mbox{ on } \Omega.

\end{array}
\right.
\end{eqnarray}

The augmented Lagrangian is  defined in $  \mathbb{V}_{\bf \epsilon}= {\bf V} \times [-\epsilon,\epsilon] \times  L^2(\Gamma)\times L^2$ by:
%
%

\begin{eqnarray*}
\label{lagrang}
\begin{array}{ccc}
\displaystyle \mathds L_\rho({\bf v},\phi,\Psi_\tau;\lambda,\mu)&=&\displaystyle J({\bf v})    + j(\Psi_\tau)+\mathds{1}_{[-\epsilon,\epsilon]}(\phi) +\left( \lambda, {\bf v}_\tau-\Psi_\tau\right)_{L^2(\Gamma)} +\\&&+\left(\mu,\nabla \cdot {\bf v}-\phi \right)+\frac{\rho}{2}\Vert \nabla \cdot {\bf v}-\phi\Vert^2+\frac{\rho}{2}\Vert {\bf v}_\tau-\Psi_\tau\Vert_{L^2(\Gamma)}^2
\end{array}
\end{eqnarray*}
%
where $\rho>0$ is the penalty parameter.  Then the problem \eqref{cst_minimi2} is equivalent to the following saddle point problem:

\begin{eqnarray}
\label{pt_selle_aug}
\left\|
\begin{array}{l}
\mbox{Find }({\bf u},\Lambda, \Phi_\tau;\lambda,\mu) \in  \mathbb{V}_{\bf \epsilon} \mbox{ such that:}\\\\
 \mathds L_\rho({\bf u},\Lambda, \Phi_\tau;\delta,\beta) \leq  \mathds L_\rho({\bf u},\Lambda, \Phi_\tau;\lambda,\mu) \leq  \mathds L_\rho({\bf v},\phi, \Psi_\tau;\lambda,\mu) ,\quad\forall ({\bf v},\phi, \Psi_\tau;\lambda)\in  \mathbb{V}_{\bf \epsilon}.
\end{array}
\right.
\end{eqnarray}

In order to obtain a saddle point that satisfy \eqref{pt_selle_aug},  we use an alternating direction method of multipliers (AMM), see \cite{fortin_glow,m} for example for more details on this variant of ADMM. This leads to the following subproblems.
Starting with an initial guess $({\bf u}^{k},{\Lambda}^k, \Phi^k_\tau;\lambda^k,\mu^k)$:
\\
compute:
\begin{eqnarray}
{\bf u}^{k+1}&\in & arg\min_{{\bf v}}   \mathds L_\rho({\bf v},\Lambda^k, \Phi^k_\tau; \lambda^k,\mu^k) \label{subp1}\\
\Lambda^{k+1}&\in & arg\min_{\phi}   \mathds L_\rho({\bf u}^{k+1},\phi, \Phi^k_\tau;\lambda^k,\mu^k) \label{subp2}\\
\Phi^{k+1}_\tau&\in & arg\min_{\Psi_\tau}   \mathds L_\rho({\bf u}^{k+1},\Lambda^{k+1}, \Psi_\tau;\lambda^k,\mu^k) \label{subp3}
\end{eqnarray}
\\
and update the Lagrange multipliers:

\begin{eqnarray}
\lambda^{k+1}&=&\lambda^{k}+\rho ({\bf u}^{k+1}_\tau-\Phi^{k+1}_\tau),\\
\mu^{k+1}&=&\mu^{k}+\rho (\nabla\cdot{\bf u}^{k+1}-\Lambda^{k+1}),
\end{eqnarray}
where $\gamma>0$ is a regularisation parameter.

To compute the solutions to subproblems \eqref{subp1}-\eqref{subp3},  we will consider  three maps.
The first considered map that is associated to the velocity is given as follows: 
\begin{eqnarray*}
V({\bf v})&:=& \mathds L_\rho({\bf v},\Lambda^k, \Phi^k_\tau;\lambda^k,\mu^k)\\&=& \displaystyle J({\bf v}) +\left( \nabla \cdot {\bf v}, \mu^k-\rho\Lambda^k\right)+\left( \lambda^k-\rho \Phi^k_\tau, {\bf v}_\tau\right)_{L^2(\Gamma)} +\frac{\rho}{2}\Vert \nabla \cdot {\bf v}\Vert^2+\frac{\rho}{2}\Vert {\bf v}_\tau\Vert_{L^2(\Gamma)}^2 +cst.
\end{eqnarray*}
where $cst$ is a constant.

The solution to the subproblems \eqref{subp1}  is computed by using the Euler-Lagrange equation as follows:
\begin{eqnarray}
\partial_u V({\bf u}){\bf v}=0.
\end{eqnarray}
where we denote by $\partial_u$ the directional derivative with respect to $u$.  This equation is equivalent to the following variational equality:
\begin{equation}
\begin{array}{c}
\displaystyle 0=\mathcal{A}({\bf u},{\bf v})-F_{p^n}({\bf v}) +\left( \nabla \cdot {\bf v}, \mu^k-\rho\Lambda^k\right)+\left( \lambda^k-\rho \Phi^k_\tau, {\bf v}_\tau\right)_{L^2(\Gamma)} +\rho\left(\nabla \cdot {\bf u}, \nabla \cdot {\bf v}\right)+\rho\left({\bf u}_\tau,{\bf v}_\tau\right)_{L^2(\Gamma)} \quad
\end{array}
\end{equation}
which is equivalent to

\begin{equation}
\begin{array}{c}
\displaystyle A({\bf u},{\bf v})+\rho\left({\bf u}_\tau,{\bf v}_\tau\right)_{L^2(\Gamma)}=F_{p^n}({\bf v})+\left( \nabla \cdot {\bf v}, \mu^k-\rho\Lambda^k\right) -\left( \lambda^k-\rho \Phi^k_\tau, {\bf v}_\tau\right)_{L^2(\Gamma)}  \quad
\end{array}
\end{equation}
where we have denoted  $A({\bf u},{\bf v})=\mathcal{A}({\bf u},{\bf v}) +\rho\left(\nabla \cdot {\bf u}, \nabla \cdot {\bf v}\right).$

The second map of the divergence is defined by:

\begin{eqnarray*}
 H(\phi)&:=& \mathds L_\rho({\bf u}^{k+1},\phi, \Phi^k_\tau;\lambda^k,\mu^k)\\&=&\frac{\rho}{2}\Vert \phi\Vert^2_0 -\left(\mu^k+\rho\nabla\cdot{\bf u}^k ,\phi\right)_{} +cst,
\end{eqnarray*}

The solution to the subproblem \eqref{subp2} can be  computed by the Euler-Lagrange inequation as follows:

\begin{eqnarray}\mbox{ Find }\Lambda\mbox{ that verify}:
\qquad(\partial_\Lambda H(\Lambda),\phi-\Lambda)\geq 0,\quad \forall \phi.
\end{eqnarray}
This gives the following variational inequality:
\begin{eqnarray}
-\left( \mu^k+\rho\nabla \cdot {\bf u}^{k+1},\phi-\Lambda\right)_{} +\rho\left(\Lambda, \phi-\Lambda\right)\geq0. \label{pro}
\end{eqnarray}
Let us put $\omega^k=\mu^k+\rho\nabla \cdot {\bf u}^{k+1}$, then \eqref{pro} can be written as follows:

\vspace{1cm}
\noindent
Find $\Lambda\in [-\epsilon,\epsilon]$  such that:
\begin{equation}
\begin{array}{cc}
\rho\left(\Lambda, \phi-\Lambda\right)\geq \left(\omega^k,\phi-\Lambda\right)_{},\forall \phi\in [-\epsilon,\epsilon]\label{pro1}
\end{array}
\end{equation}

The solution to the above variational problem is provided by the orthogonal projection of $\omega^k$ on the ball with radius $\aleph\epsilon>0$:

\begin{equation}
\Lambda^{k+1}=\dfrac{\,\omega^k}{\max(\rho\epsilon, \Vert \omega^k\Vert)}.
\end{equation}

Finally,  the third map that correspond to the friction part is given by:
\begin{eqnarray*}
 F(\Psi_\tau)&:=&  \mathds L_\rho({\bf u}^{k+1},\Lambda^{k+1}, \Psi_\tau;\lambda^k,\mu^k)\\&=& j(\Psi_\tau) -\left( \lambda^k+\rho{\bf u}_\tau^{k+1},\Psi_\tau\right)_{L^2(\Gamma)}+\rho\Vert \Psi_\tau\Vert^2_{L^2(\Gamma)}+cst.
\end{eqnarray*}

%

It is clear that the map $F$ is not Gâteaux differentiable, in  order to compute the solution to the subproblem \eqref{subp3}, we employ Fenchel duality \cite{boyd2011distributed,ekel74} and the solution is as follows:

\begin{equation}
\begin{aligned} \Theta^k & =\left\|\lambda^k+\rho {\bf u}_\tau^{k+1}\right\|_{L^2\left(\Gamma\right)}, \\ \Phi^{k+1}_\tau & =\left\{\begin{array}{cll}\frac{\Theta^k-\xi}{\rho \Theta^k}\left(\lambda^k+\rho {\bf u}_\tau^{k+1}\right), & \text { if } & \Theta^k>\xi, \\ \\0, & \text { if } & \Theta^k \leq \xi .\end{array}\right.\end{aligned}
\end{equation}


\section{Numerics}

\subsection{Algorithms and finite dimensional problem}
In this section, we discuss how our approach can address the locking phenomena.  
In the  algorithm \ref{alg1},  we summarise the solutions to  quasi-variational inequality.

\RestyleAlgo{ruled} 
\begin{algorithm}[H]
  \caption{Nearly incompressible Stokes velocity (NISV).}\label{alg1}
  \begin{enumerate}
    \item
    {\bf Initialize} $({\bf u}^{0},p^0,\Lambda^0,\Phi^0_\tau;\lambda^0,\mu^0)$.

    \item
    \While{$\Vert {\bf u}^{k+1}-{\bf u}^{k}\Vert \ge tol $,}{

    \begin{enumerate}
      \item Compute the velocity ${\bf u}^{k+1}$:
      \begin{equation}\label{velocity}
\begin{array}{c}
\displaystyle A({\bf u}^{k+1},{\bf v})+\rho\left( {\bf u}_\tau^{k+1},{\bf v}_\tau\right)_{L^2(\Gamma)}=\left( \nabla( \mu^k-\rho\Lambda^k),  {\bf v}\right)+F_{p^n}({\bf v}) +\left(\rho \Phi^k_\tau- \lambda^k, {\bf v}_\tau\right)_{L^2(\Gamma)}.  \quad
\end{array}
\end{equation}
      \item Compute the pressure $\Lambda^{k+1}$:
      
     \begin{equation}
\Lambda^{k+1}=\dfrac{\omega^{k}}{\max(\epsilon\rho, \Vert \omega^{k}\Vert)},\qquad
\end{equation}
where $\omega^k=\mu^k+\rho\nabla \cdot {\bf u}^k.$
      \item
      Compute $\Phi^{k+1}_\tau$:
      \begin{equation}\label{friction}
\begin{aligned} \Theta^{k} & =\left\|\lambda^{k}+\rho {\bf u}_\tau^{k+1}\right\|_{L^2\left(\Gamma\right)}, \\ \Phi^{k+1}_\tau & =\left\{\begin{array}{cll}\frac{\Theta^{k}-\xi}{\rho \Theta^{k}}\left(\lambda^{k}+\rho {\bf u}_\tau^{k}\right), & \text { if } & \Theta^{k}>\xi, \\\\ 0, & \text { if } & \Theta^{k} \leq \xi.\end{array}\right.\end{aligned}
\end{equation}

      \item Updates:
      \begin{eqnarray}
\lambda^{k+1}&=&\lambda^{k}+\rho ({\bf u}^{k+1}_\tau-\Phi^{k+1}_\tau),\label{multp}\\
\mu^{k+1}&=&\mu^{k}+\rho (\nabla \cdot {\bf u}^{k+1}-\Lambda^{k+1}).
\end{eqnarray}
    \end{enumerate}}

  \item
  {\bf Output}. 
  \end{enumerate}
\end{algorithm}

The solution to the whole problem \eqref{itersch} by fixed point strategy is given by the  algorithm (NIS). 

\begin{algorithm}[H]
  \caption{Nearly incompressible Stokes problem (NISP).}
  \begin{enumerate}
    \item
    {\bf Initialize} $p^0$.

    \item
    \While{$\Vert {\bf u}^{n+1}-{\bf u}^{n}\Vert_1 +\Vert p^{n+1}-p^{n}\Vert_0\ge tol $,}{
    \begin{enumerate}
      \item Compute ${\bf u}^{n+1}$ with the algorithm \ref{alg1}.

      \item Update:
      \begin{eqnarray}
p^{n+1}&=&p^{n}+\Pi_{L^2}(\nabla \cdot {\bf u}^{n+1}).
\end{eqnarray}
    \end{enumerate}}

  \item
  {\bf Output}. 
  \end{enumerate}
\end{algorithm}

In both algorithms, $tol$ is a prescribed tolerance. Another stopping criterion that  implies that are used in the algorithm (NISP) is: 
$$\Vert \nabla\cdot {\bf u}^{n+1}\Vert_0\geq tol.$$

In the  algorithm \ref{alg1}, let us investigate  the variational equation \eqref{velocity} together with the friction formula \eqref{friction} and Lagrange multiplier update \eqref{multp}.  By using the Riesz' representation theorem, it is possible to write:
 $$\left( f^k_n, {\bf v} \right)= \left(  \nabla (\mu^k-\rho\Lambda^k),{\bf v}\right)+F_{p^n}({\bf v}). $$

This enable us to write  
      \begin{equation}\label{velocity1}
\begin{array}{c}
\displaystyle A({\bf u}^{k+1},{\bf v})+\rho\left( {\bf u}_\tau^{k+1},{\bf v}_\tau\right)_{L^2(\Gamma_f)}=\left( f^k_n, {\bf v} \right) +\left( \rho \Phi^k_\tau-\lambda^k, {\bf v}_\tau\right)_{L^2(\Gamma_f)}, \quad
\end{array}
\end{equation}
and 
     \begin{equation}\label{friction1}
\begin{aligned} \Theta^{k} & =\left\|\lambda^{k}+\rho {\bf u}_\tau^{k+1}\right\|_{L^2\left(\Gamma\right)}, \\ \Phi^{k+1}_\tau & =\left\{\begin{array}{cll}\frac{\Theta^{k}-\xi}{\rho \Theta^{k}}\left(\lambda^{k}+\rho {\bf u}_\tau^{k+1}\right), & \text { if } & \Theta^{k}>\xi, \\ \\0, & \text { if } & \Theta^{k} \leq \xi, \end{array}\right.\end{aligned}
\end{equation}
with the 
updates:
      \begin{eqnarray}
\lambda^{k+1}&=&\lambda^{k}+\rho ({\bf u}^{k+1}_\tau-\Phi^{k+1}_\tau).\label{multp1}
\end{eqnarray}

The problem \eqref{velocity1}-\eqref{multp1} is  the solutions  provided by the ADMM applied  to the following partial differential equation:
\begin{eqnarray}
\label{elas_tres_cont}
\left\|
\begin{array}{cc}
-\nabla\cdot(\sigma( {\bf u}))  	=f^k_n		& \mbox{in }\Omega,\\
{\bf u}						=\textbf{0}			& \mbox{on }\Gamma_d\\
{\bf u}_{\vec{n}}=0\text{ and }
\begin{cases}|
\sigma_\tau|< \xi \Rightarrow  \quad {\bf u}_\tau	=\textbf{0},			\\
|\sigma_\tau|= \xi \Rightarrow \exists \kappa>0, \, {\bf u}_\tau	=-\kappa \sigma_\tau,	
\end{cases}

		&\mbox{on }\Gamma_f,
\end{array}
\right.
\end{eqnarray}
where $\sigma_{i j}=2 \mu \varepsilon_{i j}+\rho \delta_{i j} \varepsilon_{k k}$ is the stress tensor and $\mu=\frac{\nu}{2}$ and $\rho$ are considered here  like the Lamé coefficients.  

It is clear that the problem \eqref{elas_tres_cont} is the partial differential equation modelling  the frictional bilateral contact problem  arising in linear elasticity problem  (problem under small deformations hypothesis,  see \cite{essoufi2017alternating}) and with the volume forces $f^k_n$.  This basically ensure  that the velocity computed  at each iteration of ADMM \eqref{itersch},  is the solution to the elastic frictional contact problem \eqref{elas_tres_cont}.

In this new way, the incompressible Stokes problem with Tresca boundary condition have been approached and turned into an elastic problem without locking phenomena as it is well known that the ADMM  is convergent for finite penalty parameter $\rho>0$ (see \cite{boyd2011distributed}). Then, the locking effect will not be produced despite that  $\dfrac{1}{\epsilon}\longrightarrow +\infty$ when the velocity of the  incompressible Stokes problem is computed.  Moreover,  no inf-sup condition is needed to get the existence and uniqueness of the solution to problem in the continuous case,  so discrete inf-sup will be not required when the problem will be approximated.  All of the above arguments clearly show  that this new approach wisely made it possible to use the standard finite element approximation to approximate the problem.


\section*{Conflict of interest}

The author declare that there are no Conflict of interest.
\section*{Data availability}

This manuscript has no associated data.

\bibliographystyle{abbrv}
\bibliography{stokes__tresca__itersch_1}
\end{document}